\documentclass[a4paper]{article}

\usepackage{enumitem}
\usepackage{stmaryrd}
\usepackage{amsmath}
\usepackage{amssymb}
\usepackage{amsthm}
\usepackage[english]{babel}
\newtheorem{theorem}{Theorem}

\newtheorem{proposition}{Proposition}

\newtheorem{corollary}{Corollary}

\begin{document}

\title {Multivariate estimates for the concentration functions of weighted sums of independent identically distributed random variables}

\author{Yu.S. Eliseeva$^1$}
\begin{abstract}{In this article we formulate and prove multidimensional generalizations of the results of [6]. Let $X,X_1,\ldots,X_n$ be independent identically distributed random variables.  The paper deals with the question about the behavior of the concentration function 
of the random variable $\sum\limits_{k=1}^{n}X_k a_k$ according to 
the arithmetic structure of vectors~$a_k$.
Recently, the interest to this question has increased significantly due to the study of distributions of eigenvalues of random matrices.
In this paper we formulate 
and prove some refinements of the results of [3] and [22].}
\end{abstract}
\maketitle

\footnotetext[1]{Research supported by grants RFBR 10-01-00242, SPbU 6.38.672.2013.}

\section{Introduction}

Let $X,X_1,\ldots,X_n$ be independent identically distributed (i.i.d.) random variables. 
The concentration function of a $\mathbf{R}^d$--valued vector $Y$ with the distribution $F=\mathcal L(Y)$ is defined by the equality 
\begin{equation}
Q(F,\lambda)=\sup_{x\in\mathbf{R}^d}\mathbf{P}(Y\in x+ \lambda B), \quad \lambda>0, \nonumber
\end{equation}
where $B=\{x\in\mathbf{R}^n:\|x\|\leq 1\}$.
Let $a=(a_1,\ldots,a_n)\neq 0$, where $a_k=(a_{k1},\ldots,a_{kd})\in \mathbf{R}^d$, $k=1,\ldots, n$.
This paper deals with the question about the behavior of the concentration function of the sum $S_a=\sum\limits_{k=1}^{n} X_k a_k$ according to the arithmetic structure of vectors~$a_k$. This question is called the Littlewood--Offord problem. It is considered in [1]--[7].
The classical one--dimensional results were obtained by Littlewood and Offord [8], and Erd\"os [9] for i.i.d. $X_k$ taking values $\pm1$ with probabilities $1/2$, and integer coefficients $a_k\neq0$. In this case the concentration function is of order $O(n^{-1/2})$ (a similar estimate holds for multidimensional Littlewood--Offord problem, see [10]). However, if we assume that all $a_k$ are different, then the estimate can be significantly improved up to the order $O(n^{-3/2})$ (see [11], [12]).
 Recently, the behavior of the concentration function of weighted sums $S_a$ was actively investigated due to the study of distributions of eigenvalues of random matrices.
 
 In the sequel, let $F_a$ be the distribution of sum
$S_a=\sum\limits_{k=1}^{n} X_k a_k$, and $G$ be the distribution of a symmetrized random variable 
$\widetilde{X}=X_1-X_2$. 
Let 
\begin{equation} \label{0}M(\tau)=\tau^{-2}\int\limits_{|x|\leq\tau}x^2
\,G\{dx\}+\int\limits_{|x|>\tau}G\{dx\}=\mathbf{E}
\min\big\{{\widetilde{X}^2}/{\tau^2},1\big\}, \quad \!\tau>0.
\end{equation} 
 The symbol $c$ will be used for absolute positive constants. Note that $c$ can be different in different (or even in the same) formulas. 
We write $A\ll B$ if $|A|\leq c B$ and $B>0$. Similarly, $A\ll _{d} B$ if $|A|\leq c^d B$ and $B>0$. Note that $\ll_d$ allows constants to be exponential with respect to $d$. Also we write $A\asymp
B$ if $A\ll B$ and $B\ll A$. Similarly, $A\asymp_{d}
B$ if $A\ll_{d} B$ and $B\ll_{d} A$. For~$x=(x_1,\dots,x_d )\in\mathbf{R}^d$
we denote 
 $\|x\|^2= x_1^2+\dots +x_d^2$ and $\|x\|_\infty=\max\limits_j |x_j|$. The inner product in $\mathbf{R}^d$ is denoted $\left\langle \cdot,\cdot\right\rangle$. The product of a vector $t=(t_1,\ldots,t_d)\in\mathbf{R}^d$ and a multivector $a$  is denoted $t\cdot a=(\left\langle \, t,a_1\right\rangle,\ldots,\left\langle \, t,a_n\right\rangle)\in \mathbf{R}^n$.

Simplest properties of one-dimensional concentration functions are well studied (see, for instance, [13]--[15]). 
It is well-known that $Q(F,\mu)\ll_d (1+ \mu/\lambda)^d\,Q(F,\lambda)$ for all $\mu,\lambda>0$.
Hence, \begin{equation}\label{8a}
Q(F,c\lambda)\asymp_{d}\,Q(F,\lambda),
\end{equation}
and\begin{equation}\label{8j}
 \hbox{if }Q(F,\lambda)\ll K,\hbox{ then }Q(F,\mu)\ll_d K(1+(\mu/\lambda))^d.
\end{equation}

Recall that for any one-dimensional distribution $F$ the classical Ess\'een and Kolmogorov--Rogozin inequalities hold [16] (see as well [14] and [15]). One can find their multidimensional analogs in [17]--[20].
 
 For the random vector $Y$ with the distribution $F=\mathcal L(Y)$ in $\mathbf{R}^d$ the Ess\'een inequality holds (see Lemma 3.4, [3])
\begin{equation} \label{1n} 
Q(F,\sqrt{d}) \ll_{d} \int\limits_{B(\sqrt{d})}|\widehat{F}(t)| \,dt,
\end{equation}
where $\widehat{F}(t)=\mathbf{E}\exp(i\left\langle \, t,Y \right\rangle)$ is the characteristic function of $Y$. Let $\int\limits_{\mathbf{R}^d}|\widehat{F}(u)|\,du<\infty$ (otherwise, we can achieve this by applying smoothing), and we assume additionally that the distribution $F$ is symmetric and  $\widehat{F}(t) \geq 0$ for all~$t\in\mathbf R$, then applying relation \eqref{1n} to the measure $\cfrac{\widehat{F}(t)\,dt}{\int\limits_{\mathbf{R}^d}\widehat{F}(u)\,du}$, we obtain 
\begin{equation} \label{1a}Q(F,\sqrt{d})\gg_{d}
\int\limits_{B(\sqrt{d})} {\widehat{F}(t) \,dt}.
\end{equation}
One can find estimates of this type but with other dependence on the dimension $d$ in [21].
Thereby,
\begin{equation} \label{1b}
Q(F,\sqrt{d})\asymp_{d}\int\limits_{B(\sqrt{d})} {\widehat{F}(t) \,dt}.
\end{equation} 
 Use of relation \eqref{1b} will allow us to simplify our arguments in comparison with those in [3] and [22].
 
 Recall a multidimensional generalization of the Kolmogorov--Rogozin inequality. 

\begin{proposition}\label{thKR} Let $Y_1,\ldots,Y_n$ be independent random
vectors with the distributions $F_k=\mathcal{L}(Y_k)$.
Let $\lambda_1,\ldots,\lambda_n$ be positive numbers,
 $\lambda_k \leq \lambda$ $(k=1,\ldots,n)$. Then
\begin{equation} \label{2}
Q\Big(\mathcal{L}\big(\sum_{k=1}^{n}Y_k\big),\lambda\Big)\ll  \lambda\,
\Big(\sum_{k=1}^{n}\lambda_k^2\,\big(1-Q(\widetilde{F}_k,\lambda_k)\big)\Big)^{-1/2},
\end{equation}
where $\widetilde{F}_k$ are the distributions of corresponding symmetrized random vectors.
\end{proposition}

Siegel [18] has improved the statement of Proposition 1. He has shown that the following result holds.
\begin{proposition}\label{thE} Under the conditions of Proposition~$\ref{thKR}$ we have
\begin{equation}  \label{3}
Q\Big(\mathcal{L}\big(\sum_{k=1}^{n}Y_k\big),\lambda\Big)  \ll  \lambda \,\Big(\sum_{k=1}^{n} \lambda_k ^2 M_k(\lambda_k)\Big)^{-1/2}.
\end{equation}
\end{proposition}

One can find one-dimensional versions and refinements of these results in [13]--[15], [23]--[28].
Note that the constants in \eqref{2} and \eqref{3} are not depending on the dimension $d$. However, there exist estimates of the  Kolmogorov--Rogozin type with constants depending on $d$ (see, for example, [23], [29]).

 \medskip

The Littlewood--Offord problem was considered in [1]--[7], [22]. In this paper we formulate 
and prove multidimensional generalizations of the results of [6]. They are also the refinements of the results of [3] and [22].

 Now we formulate the results of [3] and [22] in common notation.

Friedland and Sodin [22] have simplified the arguments of Rudelson and Vershynin [2] and obtained the following result.

\begin{proposition}\label{thFS} {Let $X, X_1,\ldots,X_n$ be i.i.d. random variables and $Q(\mathcal{L}(X),1)\leq1-p$, where $p>0$, and let $a_1,\dots,a_n\in \mathbf{R}^d$. If, for any $0<D<d$ and $\alpha>0$, 
\begin{multline} \label{3b}
\sum_{k=1}^{n}(\left\langle t,a_k\right\rangle-m_k)^2\geq\alpha^2\  \hbox{for all $m_1,\dots,m_n\in \mathbf Z$, $t\in\mathbf{R}^d$}
\hbox{ such that} \\ \max\limits_{k}|\left\langle t,a_k\right\rangle|\geq1/2,\, \|t\|\leq D,
\end{multline}
 then
\begin{equation} \label{3a}
Q(F_a, d/D)\ll_{d}\exp(-cp\alpha^2)+\Big(\cfrac{\sqrt{d}}{\sqrt{p}D }\Big)^{d} \big(\det\mathbb{N}\big)^{-1/2},
\end{equation}
where 
\begin{equation} \label{33a}
 \begin{array}{rcl}
 \mathbb{N} = \sum\limits_{k=1}^{n}\mathbb{N}_k, \quad \mathbb{N}_k =  \left(\begin{array}{cccccccc}
a^{2}_{k1} & a_{k1} a_{k2} & \ldots & a_{k1} a_{kd}\\
 a_{k2} a_{k1} & a^{2}_{k2}  & \ldots & a_{k2} a_{kd}\\
 \ldots & \ldots & \ldots & \ldots\\
 \medskip
 a_{kd} a_{k1} & \ldots & \ldots & a^{2}_{kd} \end{array}\right),\\
 \\
 a_k = (a_{k1}, \ldots, a_{kd}),\quad k=1,\ldots,n.
 \end{array}
 \end{equation}}\end{proposition}

Note that the statement of Proposition \ref{thFS} in [22] was formulated and proved  in a weakened form.
There was ~$p^2$ instead of $p$ in the right-hand side of inequality~\eqref{3a}. However, the possibility to replace $p^2$ by ${p}$ was noted for example in~[3] (see Proposition \ref{thRV}). It follows easily from the elementary properties of the concentration function.

Moreover, in [22] it was assumed that  $0<D<d$. 
Furthermore there was $Q(F_a,1)$ instead of $Q(F_a, d/D)$ in the left-hand side of inequality \eqref{3a}. Since $d/D>1$ for $0<D<d$, the value $Q(F_a,1)$ can be essentially less than $Q(F_a,d/D)$. But if the authors [22] considered their result for
$D=d$, they would be deduced from it the inequality for any $D>0$ and with $Q(F_a,
d/D)$ instead of $Q(F_a,1)$ as simple as below we shall deduce Corollary
\ref{c1} from Theorem~\ref{th1}.

Note that for $\max\limits_k|\left\langle \, t,a_k\right\rangle|\leq 1/2$ we have
\begin{equation}\label{4s}
\big(\hbox{dist}(t\cdot a,\mathbf{Z}^n)\big)^2=\sum_{k=1}^{n} \min_{m_k \in \mathbf{Z}} (\left\langle t,a_k\right\rangle - m_k)^2 =
\sum_{k=1}^{n}\left\langle t,a_k\right\rangle^2, \end{equation} where
$$\hbox{dist}(t\cdot a,\mathbf{Z}^n)= \min_{m \in \mathbf{Z}^n}\|\,t\cdot a - m\|.$$
Thus, the assumption $\max\limits_{k}|\left\langle t,a_k\right\rangle|\geq1/2$ in condition~\eqref{3b} is natural.
\bigskip

Let us formulate now the multidimensional Theorem 3.3 from
 [3] in the same notation.

\begin{proposition}\label{thRV} {Let $X, X_1,\ldots,X_n$ be i.i.d. random variables with mean zero and  $Q(\mathcal{L}(X),1)\leq 1-p$, where $p>0$. 
Let $a=(a_1,\dots,a_n)$, $a_k \in \mathbf{R}^d$ such that 
$\sum\limits_{k=1}^{n}{\left\langle t, a_k\right\rangle}^2 \geq \|t\|^2$ for any $t \in \mathbf{R}^d$. 
Let $\alpha, D>0$; $\gamma\in(0,1)$, and  
\begin{equation} \label{4d}\Big(\sum_{k=1}^{n}(\left\langle t,a_k\right\rangle - m_k)^2\Big)^{1/2}\geq \min\{\gamma\|t\cdot a\|,\alpha\} \hbox{ for all }\  m_1,\ldots,m_n\in\mathbf{Z} \hbox{ and } \|t\|\leq D.\end{equation} 
Then
\begin{equation} \label{4a}Q\Big(F_a,\cfrac{d}{D}\Big) \ll_{d} \Big(\cfrac{\sqrt{d}}{\gamma D\,\sqrt{p} \, }\Big)^d +\exp(-2\,p\,\alpha^2).\end{equation} }
\end{proposition}

Note that the assumption $\mathbf{E} \,X=0$ is unnecessary in the formulation of Theorem~3.3~[3].

It is evident that if
 \begin{equation} \label{4b}
0<D\le D(a)=
\inf\big\{\|t\|>0: t\in \mathbf{R}^d, \hbox{dist}(t\cdot a,\mathbf{Z}^n)\leq \min\{\gamma\|\,t\cdot a\|,\alpha\}\big\},
\end{equation} then condition \eqref{4d} holds.
Rudelson and Vershynin [3] have called the value $D(a)$ the essential least common denominator of a vector ${a\in(\mathbf{R}^d)^n}$.

 \medskip

Now we formulate one of the main results of this paper.

\begin{theorem}\label{th1} {Let $X,X_1,\ldots,X_n$ be i.i.d. random variables. Let $a=(a_1,\ldots,a_n)$, $a_k \in \mathbf{R}^d$. Assume that, for some $\alpha>0$, condition \eqref{3b} holds for $D=\sqrt{d}$, i.e., 
\begin{multline} \label{5b}
\sum_{k=1}^{n}(\left\langle t,a_k\right\rangle-m_k)^2\geq\alpha^2\ \hbox{for all $m_1,\ldots,m_n\in \mathbf Z$, $t\in\mathbf{R}^d$}
\hbox{ such that} \\ \max_{k}|\left\langle t,a_k\right\rangle|\geq1/2, \, \|t\|\leq \sqrt{d}.
\end{multline}
 Then
$$ Q(F_a, \sqrt{d})\ll_{d} \Big(\frac{1}{\sqrt{M(1)}}\Big)^d \cfrac{1}{\sqrt{\det\mathbb{N}}}+\exp\big(-c\,\alpha^2 M(1)\big),$$
where the quantity $M(1)$ is defined in \eqref{0}, and the matrix $\mathbb{N}$ is defined in \eqref{33a}.}\end{theorem}
\medskip

Hence, it is easy to deduce what follows from Theorem \ref{th1} under the conditions of Proposition~\ref{thFS}. Namely, we have

\begin{corollary}\label{c1} {Let the conditions of Theorem\/ $\ref{th1}$ be satisfied under condition \eqref{3b} instead of~\eqref{5b} for an arbitrary $D>0$.
Then
\begin{eqnarray*}
Q\Big(F_a,\cfrac{d}{D}\,\Big) &\ll_{d}& \Big(\cfrac{\sqrt{d}}{D \sqrt{M(1)}}\,\Big)^d \cfrac{1}{\sqrt{\det\mathbb{N}}} +
\exp(-c\,\alpha^2 M(1)). 
\end{eqnarray*}}
\end{corollary}
\medskip

 Note that the value $M(1)$ is essential in refining the results of [3] and [22].  It is clear that $M(1)$ can be much larger than~$p$. For example, $p$ can be equal to $0$, but
$M(1)>0$ for any non-degenerate distribution $F=\mathcal L(X)$.  Therefore Corollary
\ref{c1} is an essential improvement of Proposition~\ref{thFS}.
It is obvious that Corollary \ref{c1} is related to Proposition~\ref{thFS} in the same way as the multidimensional variant of Ess\'een's inequality \eqref{3} is related to the multidimensional variant of the Kolmogorov--Rogozin inequality~\eqref{2}.
 
The proofs of Theorem\/ $\ref{th1}$ and Corollary \ref{c1} are in some way easier than the proofs in [3] and [22], because they do not include a complicated decomposition of the integration set. This is achieved by using relation \eqref{1b} and methods of Ess\'een [29] (see the proof of Lemma~4 of Chapter~II  in
[15]).

\medskip

We reformulate Corollary \ref{c1} for the the random variables ${X_k}/{\tau}$, $\tau>0$.

\begin{corollary}\label{c2} {Let $V_{a,\tau}=\mathcal{L}\big(\sum\limits_{k=1}^{n}a_k {X_k}/{\tau}\big)$. Then under the conditions of Corollary~ $\ref{c1}$ we have
\begin{eqnarray*}\label{4p}
Q\Big(V_{a,\tau},\cfrac{d}{D}\,\Big) &=& Q\Big(F_a,\cfrac{d \,\tau}{D}\,\Big) \ll_{d} 
\exp\left(-c\,\alpha^2 M(\tau)\right) + \Big(\cfrac{\sqrt{d}}{D\,\sqrt{M(\tau)}}\,\Big)^d \cfrac{1}{\sqrt{\det\mathbb{N}}}.
\end{eqnarray*}
For $\tau=D/d$, we obtain
\begin{eqnarray*}
Q\big(F_a,1\,\big) &\ll_{d}& \Big(\cfrac{\sqrt{d}}{D\, \sqrt{M(D/d)}}\,\Big)^d \cfrac{1}{\sqrt{\det\mathbb{N}}} + \exp\left(-c\,\alpha^2 M(D/d)\right). 
\end{eqnarray*}}
\end{corollary} 

For the proof of Corollary \ref{c2} it suffices to use relation \eqref{0}.
 \medskip

Note that $\tau$ can be arbitrarily small in Corollary \ref{c2}. Applying this statement for
$\tau$  tending to zero, we obtain
\begin{eqnarray*}Q(F_a,0) &\ll_{d}& \bigg(\cfrac{\sqrt{d}}{D \, \sqrt{{\mathbf P}(\widetilde{X}\neq 0)}}\bigg)^d \cfrac{1}{\sqrt{\det\mathbb{N}}} +
\exp\big(-c\,\alpha^2 \,{\mathbf P}(\widetilde{X}\neq 0)\big).
\end{eqnarray*}
This estimate can be deduced from the results of [3] and [22] too.
\bigskip

Now we formulate the refinements of Proposition \ref{thRV}. They are analogs of Theorem~\ref{th1} and Corollaries \ref{c1} and \ref{c2}.

\begin{theorem}\label{th2}{Let $X, X_1,\ldots,X_n$ be i.i.d. random variables.
Let $a=(a_1,\dots,a_n)$, $a_k \in \mathbf{R}^d$, $\alpha>0$, $\gamma\in(0,1)$, and  
\begin{equation} \Big(\sum_{k=1}^{n}(\left\langle t,a_k\right\rangle - m_k)^2\Big)^{1/2}\geq \min\{\gamma\|t\cdot a\|,\alpha\} \hbox{ for all }\  m_1,\ldots,m_n\in\mathbf{Z} \hbox{ and } \|t\|\leq \sqrt{d}.\end{equation} 
 Then
$$Q\big(F_a,\sqrt{d}\big) \ll_{d} \Big(\cfrac{1}{\gamma \,\sqrt{M(1)}}\,\Big)^d \cfrac{1}{\sqrt{\det\mathbb{N}}}+ \exp\left(-c\,\alpha^2 M(1)\right).$$}\end{theorem}
Note that Theorem \ref{th2} yields more general result than the result of Proposition \ref{thRV}, because the condition $\sum\limits_{k=1}^{n}\left\langle t,a_k\right\rangle ^2\geq \|t\|^2$ is absent in the formulation of Theorem \ref{th2}.

\begin{corollary}\label{c3}  {Let $X, X_1,\ldots,X_n$ be i.i.d. random variables.
Let $a=(a_1,\dots,a_n)$, $a_k \in \mathbf{R}^d$, $\alpha>0$, $D>0$, $\gamma\in(0,1)$, and  
\begin{equation} \Big(\sum_{k=1}^{n}(\left\langle t,a_k\right\rangle - m_k)^2\Big)^{1/2}\geq \min\{\gamma\|t\cdot a\|,\alpha\} \hbox{ for all }\  m_1,\ldots,m_n\in\mathbf{Z} \hbox{ and } \|t\|\leq D.\end{equation}
Then
\begin{eqnarray*}
Q\Big(F_a,\cfrac{d}{D}\,\Big) &\ll_{d}&\Big(\cfrac{\sqrt{d}}{ D\, \gamma \, \sqrt{M(1)}}\,\Big)^d \cfrac{1}{\sqrt{\det\mathbb{N}}} + \exp\left(-c\,\alpha^2 M(1)\right).
\end{eqnarray*}}
\end{corollary} 
Note that if the condition $\sum\limits_{k=1}^{n}\left\langle t,a_k\right\rangle ^2\geq \|t\|^2$ is satisfied, then the factor $\cfrac{1}{\sqrt{\det\mathbb{N}}}\leq 1$. Hence, Corollary \ref{c3} yields more general result than the result of Proposition \ref{thRV}.
Now we reformulate Corollary \ref{c3} for the variables ${X_k}/{\tau}$, $\tau>0$.

\begin{corollary}\label{c4}  {Let $V_{a,\tau}=\mathcal{L}\big(\sum\limits_{k=1}^{n}a_k {X_k}/{\tau}\big)$. 
Then under the conditions of Corollary~$\ref{c3}$ we have
\begin{eqnarray*}
Q\Big(V_{a,\tau},\cfrac{d}{D}\,\Big) &=& Q\Big(F_a,\cfrac{d \,\tau}{D}\,\Big) \ll_{d} 
\Big(\cfrac{\sqrt{d}}{D\, \gamma \,  \sqrt{M(\tau)}}\,\Big)^d \cfrac{1}{\sqrt{\det\mathbb{N}}} 
 + \exp\left(-c\,\alpha^2 M(\tau)\right).
 \end{eqnarray*}
For $\tau=D/d$, we obtain
\begin{eqnarray*}
Q\big(F_a,1\,\big) &\ll_{d}& \Big(\cfrac{\sqrt{d}}{D \,\gamma \, \sqrt{M(D/d)}}\,\Big)^d \cfrac{1}{\sqrt{\det\mathbb{N}}}  + \exp\left(-c\,\alpha^2 M(D/d)\right). 
\end{eqnarray*}}
\end{corollary}

 For the proof of Corollary \ref{c4} it suffices to use relation \eqref{0}.

\bigskip

\section{Proofs}

\emph{Proof of Theorem\/ $\ref{th1}$.} We represent the distribution $G=\mathcal{L}(\widetilde{X})$
as a mixture
$G=q E +\sum\limits_{j=0}^{\infty}p_j G_j $, where $q={\mathbf
P}(\widetilde{X}=0)$, $p_j={\mathbf
P}(\widetilde{X} \in C_j)$, $j=0,1,2,\ldots$, $C_0=\{x:|x|>1\}$, 
$C_j=\{x: 2^{-j}<|x|\leq2^{-j+1}\}$, $E$ is the probability measure concentrated at zero, $G_j$ 
are probability measures defined (for $p_j>0$) by the equality $G_j\{X\}=\cfrac{1}{p_j}
\,G\{X\bigcap C_j\}$ for any Borel set~$X$. If $p_j=0$, then we can take as $G_j$ arbitrary measures.

For $z\in \mathbf{R}$, $\gamma>0$, we introduce symmetric $d$--dimensional infinitely divisible distributions
$H_{z,\gamma}$ with the characteristic functions
\begin{equation} \widehat{H}_{z,\gamma}(t)=\exp\Big(-\cfrac{\gamma}{2}\sum_{k=1}^{n}\big(1-\cos(2\,z\left\langle \, t,a_k\right\rangle)\big)\Big), \quad \!t\in\mathbf{R}^d.
\end{equation} 
It is clear that these functions are positive everywhere.

For any characteristic function $\widehat{W}(t)$ 
of a random vector $Y$, we have
$$|\widehat{W}(t)|^2 = \mathbf{E}\exp\big(i\big\langle \, t,\widetilde{Y}\big\rangle\big) = \mathbf{E}\cos\big(\big\langle \, t,\widetilde{Y}\big\rangle\big),$$
where $\widetilde{Y}$ is a corresponding symmetrized random vector. Then
\begin{equation}\label{6}|\widehat{W}(t)| \leq
\exp\Big(-\cfrac{\,1\,}{2}\,\big(1-|\widehat{W}(t)|^2\big)\Big)  =
\exp\Big(-\cfrac{\,1\,}{2}\,\mathbf{E}\,\big(1-\cos\big(\big\langle \, t,\widetilde{Y}\big\rangle\big)\big)\Big).
\end{equation}

 Using inequalities \eqref{1n} and \eqref{6}, we obtain
\begin{eqnarray*}
Q(F_a,\sqrt{d})&\ll_{d}& \int\limits_{B(\sqrt{d})}|\widehat{F_a}(t)|\,dt \\
&\ll_{d}&
\int\limits_{B(\sqrt{d})}\exp\Big(-\frac{\,1\,}{2}\,\sum_{k=1}^{n}\mathbf{E}\,\big(1-\cos(2\left\langle \,t,a_k\right\rangle
\widetilde{X})\big)\Big)\,dt=I.
\end{eqnarray*}
It is clear that
\begin{eqnarray*}
\sum_{k=1}^{n}\mathbf{E}\big(1-\cos(2\left\langle t,a_k\right\rangle 
\widetilde{X})\big)&=&\sum_{k=1}^{n}\int\limits_{-\infty}^{\infty}\big(1-\cos(2\left\langle \,t,a_k\right\rangle x)\big)\,G\{dx\}
 \\
&=&\sum_{k=1}^{n}\sum_{j=0}^{\infty}\int\limits_{-\infty}^{\infty}\big(1-\cos(2\left\langle \,t,a_k\right\rangle x)\big)\,p_j
\,G_j\{dx\}\\ &=&\sum_{j=0}^{\infty}\sum_{k=1}^{n}\int\limits_{-\infty}^{\infty}\big(1-\cos(2\left\langle \, t,a_k \right\rangle x)\big)\,p_j \,G_j\{dx\}.
\end{eqnarray*}

We denote $\beta_j=2^{-2j}p_j $,
$\beta=\sum\limits_{j=0}^{\infty}\beta_j$, $\mu_j={\beta_j}/{\beta}$,
$j=0,1,2,\ldots$. It is evident that $\sum\limits_{j=0}^{\infty}\mu_j=1$ and
${p_j}/{\mu_j}=2^{2j}\beta $ (for $p_j> 0$).

Now we estimate the value $\beta$
\begin{eqnarray*} \label{43}
\beta = \sum_{j=0}^{\infty}\beta_j\hskip-8pt &=&\hskip-8pt\sum_{j=0}^{\infty}2^{-2j} p_j  \,
= {\mathbf P}\big(|\widetilde{X}|>1\big) +
\sum_{j=1}^{\infty}2^{-2j}\,{\mathbf
P}\big(2^{-j}<|\widetilde{X}|\leq2^{-j+1}\big)  \\
&\geq&\int\limits_{|x|>1}\,G\{dx\} + \sum_{j=1}^{\infty}
\int\limits_{2^{-j}<|x|\leq2^{-j+1}}\cfrac{x^2}{4}\,G\{dx\}\\ &\geq& \cfrac{\,1\,}{4}
\int\limits_{|x|>1}\,G\{dx\} + \cfrac{\,1\,}{4} \int\limits_{|x|\leq1}x^2 \,G\{dx\} =
\cfrac{\,1\,}{4} \,M(1).
\end{eqnarray*}
Thus,
\begin{equation}\label{9}
\beta \geq \cfrac{\,1\,}{4} \, M(1).\end{equation}

 Now we proceed like in the proof of Ess\'een's Lemma [29] (see [15], Lemma~4 of Chapter II). Applying the H\"older inequality, it is easy to see that
 \begin{equation}
 I\leq \prod _{j=0}^{\infty}I_j^{\mu_j},
 \end{equation} 
 where 
 \begin{eqnarray*}
 I_j&=&\int\limits_{B(\sqrt{d})}\exp\Big(-\cfrac{p_j}{2\,\mu_j}\;\sum_{k=1}^{n}\int\limits_{-\infty}^{\infty}\big(1-\cos(2\left\langle \, t,a_k \right\rangle x)\big)\,G_j\{dx\}\Big)\,dt \\
&=& \int\limits_{B(\sqrt{d})}\exp\Big(-2^{2j-1}\beta\;\sum_{k=1}^{n}\int\limits_{C_j}\big(1-\cos(2\left\langle \, t, a_k \right\rangle 
x)\big)\,G_j\{dx\}\Big)\,dt,
\end{eqnarray*}
if $p_j > 0$, and
 $I_j=1$ for $p_j=0$. 

Applying Ess\'een's inequality for the exponential under integral (see [15], p.~49), we have
\begin{eqnarray*}
I_j&\leq&\int\limits_{B(\sqrt{d})}\int\limits_{C_j}\exp\Big(-2^{2j-1}\beta\;
\sum_{k=1}^{n}\big(1-\cos(2\left\langle \, t,a_k \right\rangle x)\big)\Big)\,G_j\{dx\}\,dt \\
&= &\int\limits_{C_j}\int\limits_{B(\sqrt{d})}\exp\Big(-2^{2j-1}\beta\;\sum_{k=1}^{n}\big(1-\cos(2\left\langle t, a_k \right\rangle x)\big)\Big)\,dt\,G_j\{dx\} \\
&\leq& \sup_{z\in C_j}\int\limits_{B(\sqrt{d})}\widehat{H}_{z,1}^{2^{2j}\beta}(t)\,dt.
\end{eqnarray*}

We estimate the function  $\widehat{H}_{\pi,1}(t)$ for
$\max\limits_k|\left\langle \,t,a_k\right\rangle|\leq 1/2$. It is clear that there exists a $c$ such that $1-\cos x \geq cx^2$
for $|x|\leq\pi$. Thus, for $\max\limits_k|\left\langle \,t,a_k\right\rangle|\leq 1/2$,
\begin{eqnarray*}
\widehat{H}_{\pi,1}(t) &\leq& \exp\Big(-\cfrac{1}{2}\,\sum_{k=1}^{n}\big(1-\cos\left(2\pi \left\langle \, t, a_k \right\rangle\right)\big)\Big) \\ 
&\leq& \exp\Big(-c\,\sum_{k=1}^{n}|\left\langle \, t, a_k \right\rangle|^2\Big)   \leq \exp\left(-c\left\langle \mathbb{N}t,t\right\rangle\right), 
\end{eqnarray*} 
where the matrix $\mathbb{N}$ is defined in \eqref{33a}.
 
It is well-known that 
\begin{eqnarray}
\int\limits_{\mathbf{R}^d}\exp\left(-c\left\langle \mathbb{N} t, t\right\rangle \right)\, dt  \ll_{d} \cfrac{1}{\sqrt{\det \mathbb{N}}}. \label{7a}
\end{eqnarray}

For $t$ such that $
\max\limits_k|\left\langle \,t,a_k\right\rangle|\geq 1/2$, $\|t\|\leq \sqrt{d}$, one can act in the same way as in [3] and [22], namely: taking into account that 
 $$1-\cos x\geq c
\min\limits_{m\in \mathbf{Z}}|\,x-2\pi m|^2,$$ we obtain
\begin{eqnarray}
\widehat{H}_{\pi,1}(t)&\leq&\exp \Big(-c \;\sum_{k=1}^{n}\min_{m_k \in
\mathbf{Z}}\big|2\pi \left\langle \, t, a_k\right\rangle -2 \pi m_k\big|^2\Big)\nonumber \\
&=&\exp\Big(-c\;\sum_{k=1}^{n}\min_{m_k \in
\mathbf{Z}}|\,\left\langle \, t, a_k\right\rangle-m_k|^2\Big)\leq\exp(-c\,\alpha^2)\label{7b}
\end{eqnarray}
for $\|t\|\leq \sqrt{d}$ and $
\max\limits_k|\left\langle \,t,a_k\right\rangle|\geq 1/2$.

Now we will use estimates \eqref{7a} and \eqref{7b} to estimate the integrals $I_j$.
At first we consider the case $j=1,2,\ldots$. Note that the characteristic functions $\widehat{H}_{z,\gamma}(t)$ satisfy the equalities
\begin{equation} \label{5}
\widehat{H}_{z,\gamma}(t)=\widehat{H}_{y,\gamma}\big({zt}/{y}\big)\quad\hbox{and}\quad
 \widehat{H}_{z,\gamma}(t)=\widehat{H}_{z,1}^{\gamma}(t).
\end{equation}

For $z\in C_j$ we have $2^{-j}<|z|\leq2^{-j+1}<\pi$. Hence, for
${\|\,t\|\leq \sqrt{d}}$ we have $\|{zt}/{\pi}\|<\sqrt{d}$. Thus, using equalities
\eqref{5} with $y=\pi$ and estimates \eqref{7a} and \eqref{7b}, we obtain for $z\in
C_j$
\begin{eqnarray*}
\sup_{z\in C_j}\int\limits_{B(\sqrt{d})}\widehat{H}_{z,1}^{2^{2j}\beta}(t)\,dt&\leq&
\int\limits_{B(\sqrt{d})}\exp(-c\,\beta \left\langle \mathbb{N} t, t\right\rangle)\,dt +
\int\limits_{B(\sqrt{d})}\exp(-2^{2j}c\,\alpha^2 \beta )\,dt  \\
 &\ll_{d}& \Big(\cfrac{1}{\sqrt{\beta}}\Big)^d  \cfrac{1}{\sqrt{\det\, \mathbb{N}}}+ \, \exp(-c\,\alpha^2\beta).
\end{eqnarray*}

Now we consider the case $j=0$. 
Equalities \eqref{5} provide, for $z>0,\,\gamma>0$,
\begin{equation}\label{8c}
Q(H_{z,\gamma},\sqrt{d})=Q\big(H_{1,\gamma},{\sqrt{d}}/{z}\big).
\end{equation}
Thus, according to relations \eqref{8a}, \eqref{1b}, \eqref{5} and \eqref{8c}, we obtain
\begin{eqnarray*}
\sup_{z\in C_0}\int\limits_{B(\sqrt{d})}\widehat{H}_{z,1}^{\beta} (t) \,dt &=&
\sup_{z\geq 1} \int\limits_{B(\sqrt{d})}\widehat{H}_{z,\beta} (t) \,dt \asymp_{d}
\sup_{z\geq 1}\; Q(H_{z,\beta},\sqrt{d})\\ &=&
\sup_{z\geq 1}\; Q\big(H_{1,\beta},{\sqrt{d}}/{z}\big)  \leq
 Q(H_{1,\beta},\sqrt{d}) \\ &\ll_{d}&
  Q\big(H_{1,\beta},{\sqrt{d}}/{\pi}\big) =
Q(H_{\pi,\beta},\sqrt{d})\\  &\asymp_{d}&
\int\limits_{B(\sqrt{d})} \widehat{H}_{\pi,\beta}(t)\, dt =
\int\limits_{B(\sqrt{d})}\widehat{H}_{\pi,1}^{\beta}(t) \,dt.
\end{eqnarray*}

Using estimates \eqref{7a},\eqref{7b} for the characteristic function
$\widehat{H}_{\pi,1}(t)$, and the relation $\hbox{Vol}(B(\sqrt{d}))\ll_d 1$, we have:
\begin{eqnarray*}
\int\limits_{B(\sqrt{d})}\widehat{H}_{\pi,1}^\beta(t)\,dt &\leq&
\int\limits_{B(\sqrt{d})}\exp(-c\beta \left\langle \mathbb{N} t, t\right\rangle)\, dt + \int\limits_{B(\sqrt{d})}\exp(-c\,\alpha^2
\beta)\,dt \\ &\ll_{d}& \Big(\cfrac{1}{\sqrt{\beta}}\Big)^d \cfrac{1}{\sqrt{\det \mathbb{N}}} + \exp(-c\,\alpha^2 \beta).
\end{eqnarray*}

We obtained the same bound for all integrals $I_j$ for $p_j\neq 0$. In view of
$\sum\limits_{j=0}^{\infty}\mu_j=1$, we have
$$I\leq\prod_{j=0}^{\infty}I_j^{\mu_j} \ll_{d} \Big(\cfrac{1}{\sqrt{\beta}}\Big)^d \cfrac{1}{\sqrt{\det\mathbb{N}}}+
\exp(-c\,\alpha^2 \beta).$$

Hence,
\begin{eqnarray*}
Q\left(F_a,\sqrt{d}\,\right) &\ll_{d}& \Big(\cfrac{1}{\sqrt{\beta}}\Big)^d \cfrac{1}{\sqrt{\det\mathbb{N}}}+ \, \exp(-c\,\alpha^2 \beta)\\& \ll_{d}&
\Big(\cfrac{1}{\sqrt{M(1)}}\Big)^d \cfrac{1}{\sqrt{\det\mathbb{N}}} + \, \exp(-c\,\alpha^2 M(1)),
\end{eqnarray*}
that was required to prove. $\square$

\medskip

Now we deduce Corollary \ref{c1} from Theorem \ref{th1}.
\medskip

\emph{Proof of Corollary $\ref{c1}$.} We denote $$b=(b_1,\dots,b_n)=\cfrac{D}{\sqrt{d}}\,a=\cfrac{D}{\sqrt{d}}\,(a_1,\ldots,a_n) \in (\mathbf{R}^d)^n.$$ 
Then the equality $Q(F_a,{d}/{D})=Q(F_b,\sqrt{d})$ holds. The conditions of Theorem \ref{th1} for the multivector~$a$ are valid for the multivector~$b$ too. Indeed,
$\sum\limits_{k=1}^{n}(\left\langle u,b_k\right\rangle-m_k)^2\geq\alpha^2$ for all $m_1,\ldots,m_n\in \mathbf Z$, $u\in\mathbf{R}^d$ such that $\|u\|\leq \sqrt{d}$ and $\max\limits_k |\left\langle \,u,b_k\right\rangle| \geq 1/2$.
This follows from condition~\eqref{3b} of Corollary~\ref{c1}, if we denote $u=\cfrac{\sqrt{d}\,t}{D}$.
It remains to apply Theorem \ref{th1} to the multivector $b$. $\square$
\bigskip

\emph{Proof of Theorem\/ $\ref{th2}$.} We will act similarly to the proof of Theorem \ref{th1}. Using the notation of Theorem
\ref{th1}, we recall that
\begin{eqnarray*}
Q\left(F_a,\sqrt{d}\,\right) &\ll_d& \prod_{j=0}^{\infty}\sup_{z \in C_j}\int\limits_{B(\sqrt{d})}\widehat{H}_{z,1}^{2^{2j}\beta }(t)\,dt \\ &\leq& \prod_{j=0}^{\infty}\sup_{z \in C_j}\int\limits_{B(\sqrt{d})}\widehat{H}_{\pi,1}^{2^{2j}\beta }\big({zt}/{\pi}\big)\,dt.
\end{eqnarray*}

The conditions of Theorem \ref{th2} imply that
\begin{eqnarray*}
\widehat{H}_{\pi,1}(t) &\leq &\exp\Big(-c\;\sum_{k=1}^{n}\min_{m_k\in\mathbf{Z}}\;
\bigl|2\pi \langle \,t, a_k\rangle - 2\pi m_k\bigr|^2\Big)\\ &\leq& \exp(-c\,\alpha^2) + \exp\big(-C\,\gamma^2\langle \mathbb{N} t, t\rangle\big)
\end{eqnarray*} 
for all  $\|t\| \leq \sqrt{d}$, where $\mathbb{N}$ is defined in \eqref{33a}.
Hence,
\begin{eqnarray*}
Q\left(F_a,\sqrt{d}\,\right) &\ll_d& \int\limits_{B(\sqrt{d})}\exp(-c\,\gamma^2\beta\, \left\langle \mathbb{N} t, t\right\rangle)\,dt + \int\limits_{B(\sqrt{d})}\exp(-c\,\alpha^2\beta)\,dt \\ &\ll_d& \Big(\cfrac{1}{\gamma\,\sqrt{\beta}}\Big)^d \cfrac{1}{\sqrt{\det\mathbb{N}}} + \exp(-c\,\alpha^2 \beta).
\end{eqnarray*}
According to \eqref{43}, $\beta \geq M(1)/4$. Then we obtain
$$Q\left(F_a,\sqrt{d}\,\right) \ll_d \Big(\cfrac{1}{\gamma\, \sqrt{M(1)}}\Big)^d \cfrac{1}{\sqrt{\det\mathbb{N}}}+ \exp(-c\,\alpha^2M(1)),$$
that was required to prove. $\square$\medskip

\emph{Proof of Corollary\/ $\ref{c3}$.} This proof is similar to the proof of Corollary \ref{c1}.
We denote $b=\cfrac{D}{\sqrt{d}}\,a \in (\mathbf{R}^d)^n$ and $u=\cfrac{\sqrt{d}\,t}{D}$. Then 
$\Big(\sum\limits_{k=1}^{n}(\left\langle u, b_k\right\rangle - m_k)^2\Big)^{1/2}\geq \min\{\gamma\|t\cdot a\|,\alpha\}$ for all  $m_1,\ldots,m_n\in\mathbf{Z}$ and $\|u\|\leq \sqrt{d}$. Thus, the conditions of Theorem \ref{th2} for the multivector~$a$ are valid for the multivector~$b$ as well. It remains to note that $Q\big(F_a,{d}/{D}\big)=Q(F_b,\sqrt{d})$ and to apply Theorem  \ref{th2} to the multivector $b$. $\square$
\bigskip
\bigskip
\bigskip
\bigskip

\begin{center} References \end{center}
{\nopagebreak}
\small

\begin{enumerate}[topsep=0pt, itemsep=-0.5ex]

\item
Nguyen H., Vu V., \textit{Optimal inverse Littlewood--Offord theorems.}
 Adv. Math. \textbf{226} (2011), 5298--5319. 

\item
Rudelson M., Vershynin R., \textit{The Littlewood--Offord problem and invertibility of random matrices.}
 Adv. Math.  \textbf{218} (2008), 600--633.

\item
Rudelson M., Vershynin R., \textit{The smallest singular value of a random rectangular matrix.}
 Comm. Pure Appl. Math. \textbf{62} (2009), 1707--1739.

\item
Tao T., Vu V., \textit{Inverse Littlewood--Offord theorems and the
condition number of random discrete matrices.}
 Ann. Math. \textbf{169} (2009), 595--632.

\item
Tao T., Vu V., \textit{From the Littlewood--Offord problem to the circular law: universality of the spectral distribution of random matrices.}
 Bull. Amer. Math. Soc. \textbf{46} (2009), 377--396.

\item
Eliseeva Yu. S., Zaitsev A. Yu., \textit{Estimates for the concentration function of weighted sums of independent random variables.}  
 Theory Probab. Appl. \textbf{57} (2012), 768--777.
 
 \item
G. Hal\'asz, \textit{Estimates for the concentration function of combinatorial number theory and probability.}
--- Periodica Mathematica Hungarica \textbf{8} (1977), 197--211. 
 
\item
Littlewood J. E., Offord A. C., \textit{On the number of real roots of a random algebraic equation.}
 Rec. Math. [Mat. Sbornik] N.S. \textbf{12} (1943), 277--286.

\item
Erd\"os P., \textit{On a lemma of Littlewood and Offord.}
 Bull. Amer. Math. Soc. \textbf{51} (1945), 898--902.

\item
Frankl P., F\"uredi Z., \textit{Solution of the Littlewood--Offord problem in high dimensions.}
 Ann. Math. \textbf{128} (1988), 259--270.

\item
Erd\"os P., \textit{Extremal problems in number theory.}
 Amer. Math. Soc., Providence, R.I. \textbf{VIII} (1965), 181--189 

\item
S\'ark\"ozy A., Szemer\'edi E., \textit{Uber ein problem von Erd\"os und Moser.}
 Acta Arithmetica \textbf{11} (1965), 205--208.

\item
Arak T. V., Zaitsev A. Yu., \textit{Uniform limit theorems for sums of independent random variavles.} 
 Proc. Steklov Inst. Math. p. 174 (1988)
\item
Hengartner W., Theodorescu R., \textit{Concentration function.}
 Academic Press, New York, (1973).

\item
Petrov V. V., \textit{Sums of independent random variables.}
 Moscow, Science, (1972).

\item
Ess\'een C.-G., \textit{On the Kolmogorov--Rogozin inequality for the concentration function.}
 Z. Wahrscheinlichkeitstheor. verw. Geb. \textbf{5} (1966), 210--216.

\item
Tao T., Vu V., \textit{Additive Combinatorics.} 
 Cambridge University Press, Cambridge \textbf{105} (2006).

\item
Siegel G., \textit{Upper bounds for tre concentration function in a Hilbert space.}
 Theory Probab. Appl. \textbf{26} (1981), 335--349.

\item
Miroshnikov A. L., \textit{Bounds for the multidimensional L\'evy concentration function.}
 Theory Probab. Appl. \textbf{34} (1989), 535--540.

\item
Anan'evskii S. M., Miroshnikov A. L., \textit{Local bounds for the L\'evy concentration function in a multidimensional or a Hilbert space.}
  Zap. Nauchn. Sem. LOMI \textbf{130} (1983), 6--10.

\item
Zaitsev A. Yu., \textit{For the multidimensional generalization of the method of triangular functions.}
 Zap. Nauchn. Sem. LOMI \textbf{158} (1987), 81--104.

\item
Friedland O., Sodin S., \textit{Bounds on the concentration function
in terms of Diophantine approximation.}
 C. R. Math. Acad. Sci. Paris \textbf{345} (2007), 513--518.

\item
Rogozin B. A., \textit{On the increase of dispersion of sums of independent random variables.}
 Theory Probab. Appl. \textbf{6} (1961), 97--99.

\item
Arak T. V., \textit{On the convergence rate in Kolmogorov's uniform limit theorem. I.}
Theory Probab. Appl. \textbf{26} (1981), no. 2, 225--245..

\item
Bretagnolle J., \textit{Sur l'in\'egalit\'e de concentration
de Doeblin--L\'evy, Rogozin--Kesten.}  In: Parametric and
semiparametric models with applications to reliability, survival
analysis, and quality of life, Stat. Ind. Technol., Boston: Birkh\"auser
(2004), 533--551.

\item
Kesten H., \textit{A sharper form of the
Doeblin--L\'evy--Kolmogorov--Rogozin inequality for concentration
functions.}
 Math. Scand. \textbf{25} (1969), 133--144.

\item
Miroshnikov A. L., Rogozin B. A., \textit{Inequalities for the concentration function.}
 Theory Probab. Appl. \textbf{25} (1980), 176--180.

\item
Nagaev S. V., Hodzhabagyan S. S., \textit{On the estimate for the concentration function of sums of independent random variables.}
 Theory Probab. Appl. \textbf{41} (1996), 655--665.

\item
Ess\'een C.-G., \textit{On the concentration function of a sum of
independent random variables.}
 Z. Wahrscheinlichkeitstheor. verw.
Geb. \textbf{9} (1968), 290--308.

\end{enumerate}
\end{document}